\documentclass[10pt]{article}
\usepackage{amsfonts}
\usepackage{mathrsfs}
\usepackage{amsmath}
\usepackage{amssymb}
\usepackage{graphicx}
\usepackage{epic}
\usepackage{amsthm,latexsym}
\renewcommand{\paragraph}{\roman{paragraph}}
\def \dist{\hbox{dist}}
\def \P{{\bf P}}
\def \C{{\bf C}}
\usepackage[]{caption2} 
\newtheorem{theorem}{\scshape \mdseries  Theorem}[section]
\newtheorem{lemma}[theorem]{\scshape \mdseries  Lemma}
\newtheorem{coro}[theorem]{\scshape \mdseries  Corollary}

\begin{document}

\title{\sf The triangle-free graphs with rank $6$\thanks{Supported by National Natural Science Foundation of China (11071002),
Program for New Century Excellent Talents in University, Key Project of Chinese Ministry of Education (210091),
Specialized Research Fund for the Doctoral Program of Higher Education (20103401110002),
Anhui Provincial Natural Science Foundation (10040606Y33),
Scientific Research Fund for Fostering Distinguished Young Scholars of Anhui University(KJJQ1001),
Academic Innovation Team of Anhui University Project (KJTD001B),
 Fund for Youth Scientific Research of Anhui University(KJQN1003).}}
\author{Long Wang, Yi-Zheng Fan\thanks{Corresponding author. E-mail address: fanyz@ahu.edu.cn (Y.-Z. Fan), wangy@ahu.edu.cn (Y. Wang), wanglongxuzhou@126.com (L. Wang)},
\ Yi Wang\\
  {\small  \it School of Mathematical Sciences, Anhui University, Hefei 230601, P. R. China}
 }
\date{}
\maketitle

\noindent {\bf Abstract:} The rank of a graph $G$ is defined to be the rank of its adjacency matrix $A(G)$.
In this paper we characterize all connected triangle-free graphs with rank $6$.

\noindent {\bf MR Subject Classifications:} 05C50


\section{Introduction}
Throughout this paper we consider simple graphs.
The {\it rank of a graph} $G=(V(G),E(G))$, denoted by $r(G)$, is  defined to be the rank of its adjacency matrix;
 and the {\it nullity of $G$}, denoted by $\eta(G)$, is defined to be the multiplicity of zero eigenvalues of its adjacency matrix.
It is easy to see that $r(G)+\eta(G)=|V(G)|$.
The graph $G$ is called {\it singular} or {\it nonsingular} if $\eta(G)>0$ or $\eta(G)=0$.

In chemistry, a conjugated hydrocarbon molecule can be modeled by its molecular graph.
 It was known in $\cite{cvet}$ or $\cite{lon}$, if the molecule represented by a graph $G$ is chemically stable, then $G$ is nonsingular.
 In 1957 Collatz and Sinogowitz $\cite{colla}$ posed the problem of characterizing all singular graphs.
 The problem is very hard and only some particular results are known (see $\cite{che2}$, $\cite{guo}$, $\cite{hu}$, $\cite{tang}$), although it  has received a lot of attention.

 We review some known results related to this topic.
 In \cite{fio} and \cite{li} the smallest rank among $n$-vertex trees in which no vertex has degree greater than a fixed value was determined,
  and the corresponding trees were constructed.
  The  singular line graphs of trees were described in \cite{gut} and \cite{sci}.
  The rank set of  bipartite graphs of fixed order and bipartite graphs of rank $4$ were determined in \cite{fan}.
It was shown in  \cite{che1} that all connected graphs with rank $2$ (respectively, $3$) are complete bipartite graphs (respectively,  complete tripartite graphs).
 In \cite{cha1} and \cite{cha2}, the authors proved that each connected graph with rank $4$ (or $5$) is obtained from one of the $8$ (or $24$) reduced graphs by multiplication of vertices.

Now a natural problem is left open: {\it To characterize graphs with rank $6$.}
Because a graph of rank $6$ can be obtained by multiplication of vertices from its reduced form, the problem is equivalent to characterize reduced graphs with rank $6$.
We also note that regular bipartite graphs with rank $6$ has been characterized in $\cite{fan}$.
In the present paper, we focus our attention on reduced connected triangle-free graphs with rank $6$, and characterize all such graphs.

\section{Preliminaries}

We first introduce two graph operations: {\it multiplication and reduction of vertices}.
Given a graph $G$ on vertices $v_{1}, v_{2},\ldots,v_{n}$.
Let $m=(m_{1},m_{2},\ldots,m_{n})$ be a list of positive integers.
Denote by $G\circ m$ the graph obtained from $G$ by replacing each vertex $v_{i}$ of $G$ with an independent set of $m_{i}$ vertices $v_{i}^{1},v_{i}^{2},\ldots,v_{i}^{m_{i}}$,
 and joining $v_{i}^{s}$ with $v_{j}^{t}$ if and only if $v_{i}$ and $v_{j}$ are adjacent in $G$.
 The resulting graph $G\circ m$ is said to be obtained from $G$ by {\it multiplication of vertices}; see \cite{cha1,cha2}.

Define a relation $\approx$ in $V(G)$ in the way that $x\approx y$ if and only if $N_G(x)=N_G(y)$, where $N_G(x)$ denotes the neighborhoods of a vertex $x$ in $G$.
Obviously, the relation is an equivalence one, and each equivalence class $\bar v=\{x:\ N_G(x)=N_G(v)\}$ is an independent set.
Now construct a new graph  $R(G)$  obtained from $G$ by taking the vertex set to be  all equivalence classes $\bar v$ and joining $\bar v$ with $\bar u$ if and only if $v, u$ are adjacent in $G$.
The graph $R(G)$ is called to be obtained from $G$ by {\it reduction of vertices}, and is a {\it reduced form } of $G$.
 One can see the above two operations are inverse to each other, and preserve the rank of graphs.

A graph is called {\it reduced} if itself is a reduced form, i.e. the neighborhoods of distinct vertices are distinct.
 Denote by $H \subseteq G$ if $H$ is a subgraph of $G$, and $H \lhd G$ if $H$ is an induced subgraph of $G$.
 For $H \subseteq G$ and $v\in V(G)$, denote by $N_{H}(v)$ the neighborhoods of $v$ in $H$.
 If two vertices $u$ and $v$ are adjacent in $G$, then we write $u\sim v$; otherwise write $u \nsim v$.
 The  distance of two vertices $u, v$ in $G$ is denoted as $\dist_{G}(u,v)$.
 For $w\in V(G)\backslash V(H)$ and $H\lhd G$, the distance between $w$ and $H$ is denoted and defined by $\dist_{G}(w,H)=\min_{v\in V(H)}\{\dist_{G}(w,v)\}$.
 For later use we now introduce some basic results.

\begin{lemma}\label{tree-nul}{\em\cite{cvet}} Let $G$ be a tree, then $r(G)=2\mu(G)$, where $\mu(G)$ is the matching number of $G$.\end{lemma}

\begin{lemma}\label{pend-nul} {\em\cite{cvet}}
 Let $G$ be a graph containing a pendant vertex, and let $H$ be the induced subgraph of G by deleting the pendant vertex and the vertices adjacent to it.
 Then $\eta(G)=\eta(H)$, or equivalently, $r(G)=r(H)+2$.
 \end{lemma}

\begin{lemma}\label{bi-nulset} {\em \cite{fan}} Let $\mathcal{B}_{n}$ be the set of bipartite graphs of order $n$.
The nullity set of $\mathcal{B}_{n}$ is $\{n-2k:k=0,1,2,\ldots, \lfloor\frac{n}{2}\rfloor$\}.\end{lemma}

 \begin{lemma}\label{neigh-nonaj}{\em \cite{wong}} Let $H$ be an induced subgraph of a reduced graph $G$ for which $r(H)=r(G)$.
 If $u,v\in V(G)\backslash V(H)$ are not adjacent in $G$, then $N_{H}(u)\neq N_{H}(v)$.
 \end{lemma}

 \begin{lemma}\label{neigh-outside}{\em \cite{wong}}
 Let $H$ be an induced subgraph of a reduced graph $G$ for which $r(H)=r(G)$, and let $v$ be a vertex not in $H$.
 Then $N_{H}(v)\neq N_{H}(u)$ for any $u\in V(H)$.
 \end{lemma}

 \begin{lemma}\label{dist} {\em \cite{wong}} Let $H$ be a proper induced subgraph of a connected graph $G$ for which $r(H)\geq r(G)-1$.
 Then $\dist_{G}(v,H)=1$ for each vertex $v\in V(G)\backslash V(H)$.
 \end{lemma}


 \begin{lemma}\label{op1}
 Let $H$ be a nonsingular induced subgraph of a reduced graph $G$ for which $r(H)=r(G)$.
  Assume $G$ has exactly two distinct vertices $u,v$ outside $H$.
  Let $\tilde{G}$ be obtained from $G$ by adding a new edge $uv$ if $uv \notin E(G)$ or deleting the edge $uv$ otherwise. Then $r(\tilde{G})=r(G)+2$.
  \end{lemma}

 \noindent{\bf Proof.}
The adjacency matrix of $G$ and $\tilde{G}$ can be written as:
 $$A(G)=\left(
          \begin{array}{ccc}
            0 & \theta & \alpha \\
            \theta & 0 & \beta \\
            \alpha^{T} & \beta^{T} & B \\
          \end{array}
        \right), \  A(\tilde{G})=\left(
                                \begin{array}{ccc}
                                  0 & 1-\theta & \alpha \\
                                  1-\theta & 0 & \beta \\
                                  \alpha^{T} & \beta^{T} & B \\
                                \end{array}
                              \right),$$
 where the first two rows of both matrices correspond to $u,v$ respectively,
 $\theta$ equals $1$ or $0$ if $uv \in E(G)$ or not, $B$ is the adjacency matrix of $H$.
  Let $Q=\left(
                                                                                 \begin{array}{ccc}
                                                                                   1 & 0 & 0 \\
                                                                                   0 & 1 & 0 \\
                                                                                   -B^{-1}\alpha^{T} & -B^{-1}\beta^{T} & I_{k} \\
                                                                                 \end{array}
                                                                             \right).$
 Then {\small
 $$Q^{T}A(G)Q=\left(
                \begin{array}{ccc}
                  -\alpha B^{-1} \alpha^T & \theta-\alpha B^{-1}\beta^{T} & 0 \\
                  \theta-\beta B^{-1}\alpha^{T} & -\beta B^{-1}\beta^{T} & 0 \\
                  0 & 0 & B \\
                \end{array}
              \right)
 ,Q^{T}A(\tilde{G})Q=\left(
     \begin{array}{ccc}
       -\alpha B^{-1} \alpha^T & 1-\theta-\alpha B^{-1}\beta^{T} & 0 \\
       1-\theta-\beta B^{-1}\alpha^{T} & -\beta B^{-1}\beta^{T} & 0 \\
       0 & 0 & B \\
     \end{array}
   \right).
 $$}
Since $r(G)=r(H)=r(B)$, the left-upper submatrix of order $2$ of $Q^{T}A(G)Q$ is zero.
So
$$Q^{T}A(\tilde{G})Q=\left(
                       \begin{array}{ccc}
                         0 & 1 & 0 \\
                         1 & 0 & 0 \\
                         0 & 0 & A(H) \\
                       \end{array}
                     \right)
,$$
whose rank is $k+2$.\hfill$\blacksquare$

\begin{coro}\label{op2}
Let $H$ be a nonsingular induced subgraph of a reduced graph $G$ for which $r(H)=r(G)$.
 Assume that $V(G)\backslash V(H)=\{v_{1},v_{2},\ldots,v_{k}\}$, $k\geq 2$.
 Let $\bar{G}$ be obtained from $G$ by adding some edges if these edges are not in $G$ or deleting some edges if they are in $G$.
Then $r(\bar{G})\geq r(G)+2$.
\end{coro}

 \noindent{\bf Proof.}  Without loss of generality, the adjacency relation of $v_1,v_2$ is different between $G$ and $\bar{G}$.
 Taking the subgraph $G'$ of $G$ and $\bar{G}'$ of $\bar{G}$ both induced by $V(H)\cup \{u,v\}$, we have
 $r(\bar{G}') \ge r(G')+2$ by Lemma \ref{op1}.
 The result now follows as $r(\bar{G}) \ge r(\bar{G}')$ and $ r(G')=r(G)$.\hfill$\blacksquare$

\section{Characterizing nonsingular triangle-free graphs of order $6$}

 Let $\P_{k}$ denote the path of order $k$, and let $\C_{k}$ denote the cycle of order $k$ (to avoid the confusion of $C_k$ used in Fig. \ref{picC}).
 The disjoint union of two graphs $G$ and $H$ is written as $G\cup H$.
 Denote by $kH$  the disjoint union of $k$ copies of $H$.
  In this section, we obtain the following main result by proving a series of lemmas below.

\begin{figure}[h!]
  \renewcommand\thefigure{\arabic{section}.\arabic{figure}}
   \centering
  \includegraphics[scale=0.7]{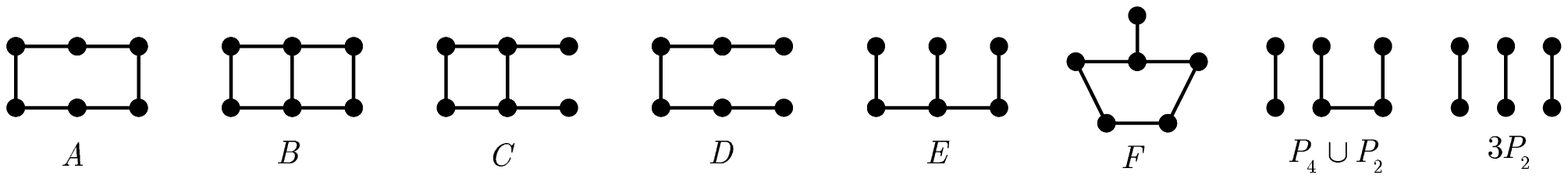}\\
  \caption{\small $8$ nonsingular graphs $A,B,C,D,E,F,3\P_{2},\P_{2}\cup \P_{4}$}\label{pic-nons}
\end{figure}

\begin{theorem} \label{nons-6}
If $G$ is a nonsingular triangle-free graph of order $6$, then $G$ is one graph in Fig. \ref{pic-nons}.
\end{theorem}

For a nonsingular bipartite graph $G$ of order $6$, if $G$ is disconnected, by Lemma \ref{bi-nulset},
the rank of each component must be $2$ or $4$.
Note that connected nonsingular bipartite graphs of order $2$ must be $\P_{2}$,
 and connected nonsingular bipartite graphs of order $4$ must be $\P_{4}$; see $\cite{cha1}$ or $\cite{che1}$ for details.
 So we have the following result immediately.

\begin{lemma} \label{nons-6-disc}
Let $G$ be a disconnected nonsingular bipartite graph of order $6$. Then $G$ is $3\P_{2}$ or $\P_{4}\cup \P_{2}$.
\end{lemma}

Denote by $g(G)$ the {\it girth} of a graph $G$ (i.e., the shortest length of cycles in $G$).
If $G$ is a nonsingular triangle-free graph of order $6$, then either $g(G)=4,5,6$, or $G$ is a tree. If $g(G)=6$, clearly $G=C_6$.

\begin{lemma}\label{nons-6-g5}
If $G$ is a connected nonsingular graph of order $6$ with $g(G)=5$, then $G=F$ in Fig. \ref{pic-nons}.
\end{lemma}

\noindent{\bf Proof.} Let $v$ be the unique vertex of $G$ outside $\C_{5}$.
If $|N_{\C_{5}}(v)|\geq 3$, then $G$ contains triangles.
If $|N_{\C_{5}}(v)|=2$, say $N_{\C_{5}}(v)=\{u_{1},u_{2}\}$, then $\dist_{\C_{5}}(u_{1},u_{2})= 2$; otherwise $G$ contains triangles.
However, in this case $G$ is not reduced and hence has rank at most $5$. Thus $|N_{\C_{5}}(v)|=1$ and the result follows. \hfill$\blacksquare$

\begin{lemma}\label{nons-6-g4}
If $G$ is a connected nonsingular graph of order $6$ with $g(G)=4$, then $G$ is the graph $B$ or $C$ in Fig. \ref{pic-nons}.
\end{lemma}

\noindent{\bf Proof.} Note that $r(C_{4})=2$, $r(G)=6$, and deleting a vertex from $G$ reduces $r(G)$ at most $2$.
There exists a graph $H \lhd G$ such that $|V(H)|=5$ and $r(H)=4$. As shown in $\cite{cha1}$, $H$ is obtained from $\C_{4}$ by attaching a pendant vertex.
For the reconstruction of $G$, one need to add a vertex together with some edges.
It is easily to check that $G$ is the graph $B$ or $C$.\hfill$\blacksquare$

\begin{lemma} \label{nons-6-tree}
If $G$ is a nonsingular tree of order $6$, then $G$ is the graph $D$ or $E$ in Fig. \ref{pic-nons}.
\end{lemma}

\noindent{\bf Proof.} By Lemma \ref{tree-nul}, $\mu(G)=3$. Thus $G \supseteq 3\P_{2}$.
Since a tree of order $6$ must contain $5$ edges, we have $|E(G\backslash 3\P_{2})|=2$.
One can easily see that $G$ is the graph $D$ or $E$.  \hfill$\blacksquare$

\section{Characterizing reduced triangle-free graphs with rank $6$}
Let $G$ be a connected reduced triangle-free graph with rank $6$.
We know that $G$ contains an induced subgraph of order $6$ and rank $6$, which are listed in Fig. \ref{pic-nons} by Theorem \ref{nons-6}.
In order to reconstruct $G$ we are left to consider how to add vertices to the graphs in Fig. \ref{pic-nons}.
The main result of this paper is as follows.

\begin{theorem}\label{rank6} Let $G$ be a connected reduced triangle-free graph.
 Then $r(G)=6$ if and only if there exist two graphs $H_{1}$ and $H_{2}$ such that
 $H_{1}\rhd G\rhd H_{2}$, where $H_1$ is one graph in Fig. \ref{pic-nons} and $H_2$ is one graph in Fig. \ref{pic-rank-6}.
\end{theorem}

\setcounter{figure}{0}
\begin{figure}[h!]
\renewcommand\thefigure{\arabic{section}.\arabic{figure}}
  \centering
  \includegraphics[scale=0.7]{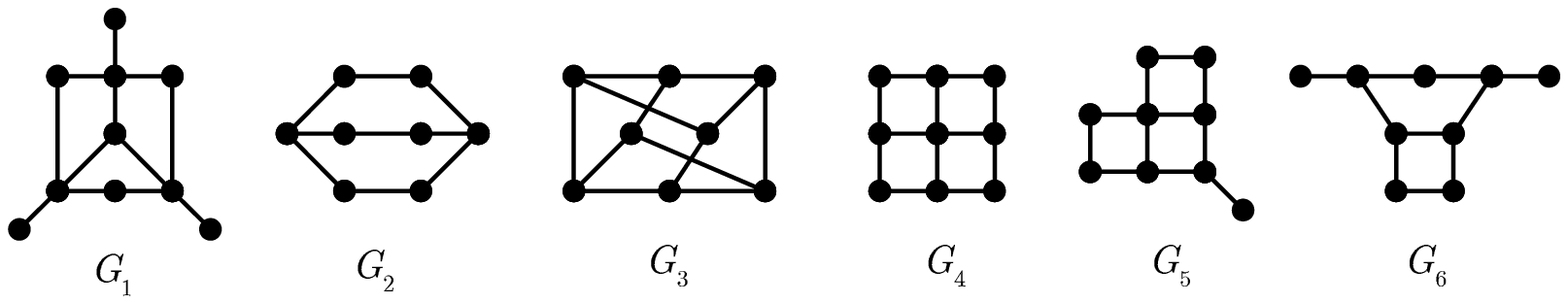}\\
  \caption{\small $6$ reduced graphs with rank $6$: $G_{1},G_{2},G_{3},G_{4},G_{5},G_{6}$}\label{pic-rank-6}
\end{figure}

It is easy to check each graph in Fig. \ref{pic-rank-6} has rank $6$.
So the sufficiency of Theorem \ref{rank6} follows.
We first discuss the necessity for the case of $G$ being bipartite, and then the case of
$G$ being non-bipartite.
For $H\lhd G$, Denote by $H(k)$ the set of vertices $v \in V(G) \backslash V(H)$ such that $|N_{H}(v)|=k$.

\subsection{Reduced bipartite graphs of rank $6$}

\begin{lemma}\label{degree}
Let $G$ be a connected reduced bipartite graph with rank $6$.\\
(1) If $G\supseteq H_{1}\supseteq \P_{6}$, where $|V(H_{1})|=6$, then $|N_{H_{1}}(v)|\leq 3$ for each $v\in V(G)\backslash V(H_{1})$;\\
(2) If $G\supseteq H_{2}\supseteq C$ of Fig. \ref{pic-nons}, where $|V(H_{2})|=6$, then $|N_{H_{2}}(v)|\leq 2$ for each $v\in V(G)\backslash V(H_{2})$.
\end{lemma}

\noindent{\bf Proof.} Suppose that $v\in V(G)\backslash V(H_{1})$ and $N_{H_{1}}(v)=\{v_{1},v_{2},\ldots,v_{k}\} \ne \emptyset$.
Then $\dist_{H_{1}}(v_{i},v_{j})$ is even; otherwise $G$ is not bipartite.
As $H_{1}\supseteq \P_{6}$, there exist at most three vertices whose pairwise distance is even in $H_{1}$, which implies that  $k\leq 3$.
The proof for the next result is similar and omitted.\hfill$\blacksquare$

\begin{lemma} \label{bi-A}
Let $G$ be a connected reduced bipartite graph with rank $6$, which contains $A$ (or $\C_6$) in Fig. \ref{pic-nons} as an induced subgraph.
Then $G$ is an induced subgraph of $G_{1}$, $G_{2}$ or $G_{3}$ in Fig. \ref{pic-rank-6}
\end{lemma}

\begin{figure}[h!]
\renewcommand\thefigure{\arabic{section}.\arabic{figure}}
  \centering
  \includegraphics[scale=0.7]{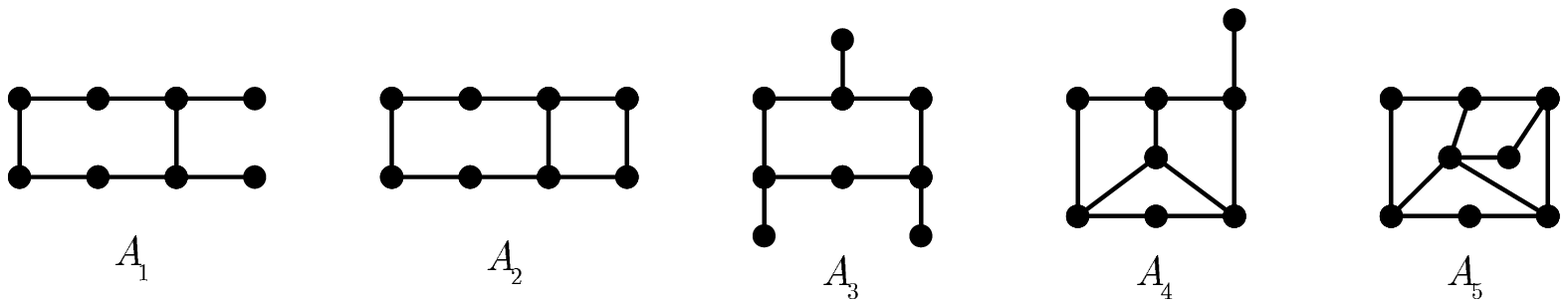}\\
  \caption{\small Illustration for Lemma \ref{bi-A}, where $r(A_{1})=r(A_{2})=r(A_{4})=r(A_{5})=8$, $r(A_{3})=6$}\label{picA}
\end{figure}

\noindent{\bf Proof.}
By Lemma \ref{dist} and Lemma \ref{degree}, $1\leq |N_{A}(v)|\leq 3$ for $v\in V(G)\backslash V(A)$.
If $A(2)\neq\emptyset$, let $v\in A(2)$ and $v_{1}$, $v_{2}$ be its two neighbors in $A$.
If $\dist_{A}(v_{1},v_{2})$ is odd, then $G$ contains an odd cycle.
If $\dist_{A}(v_{1},v_{2})=2$, then $G$ is not reduced by Lemma \ref{neigh-outside}. So $A(2)=\emptyset$. We divide the remaining proof into some cases:

{\it Case 1:} $A(3)=\emptyset$.
If $|A(1)|=1$, then $G\lhd G_{1}$ of Fig. \ref{pic-rank-6} obviously. Now we assume $|A(1)|\geq 2$.
Suppose $v_{1},v_{2},\ldots,v_{k}\in A(1)$ and $v_{1}',v_{2}',\ldots,v_{k}'$ are their unique neighbors in $A$, respectively.
Note that $v_{i}'\neq v_{j}'$ for $i\neq j$ as $G$ is reduced.
 Observing the graph $A_1,A_2$ in Fig. \ref{picA} both have rank $8$, so $\dist_{A}(v_{i}',v_{j}')\neq 1$ for $i\neq j$,
 which implies that $\dist_{A}(v_{i_{1}}',v_{j_{1}}')=2$ and $\dist_{A}(v_{i_{2}}',v_{j_{2}}')=3$ can not hold at the same time
  for some $v_{i_{1}}',v_{i_{2}}',v_{j_{1}}',v_{j_{2}}'$.
 If $\dist_{A}(v_{i}',v_{j}')=2$ for some $i\neq j$, then $|A(1)|\leq 3$.
 Since $\C_6\lhd A_3 \lhd G_1$, we have $r(A_{3})=6$, where $A_3$ is listed in Fig. \ref{picA}.
 It follows that $v_{i} \nsim v_{j}$ in $G$; otherwise $r(G) \ge 8$ by Lemma \ref{op1}.
 Hence $G\lhd A_{3}\lhd G_{1}$.
 If $\dist_{A}(u_i',v_j')=3$, then obviously $|A(1)|= 2$.
 Since $r(G_{2})=6$, $v_{i} \sim v_{j}$ in $G$ also by Lemma \ref{op1}, which implies that $G=G_{2}$, where $G_2$ is listed in Fig. \ref{pic-rank-6}.

{\it Case 2:} $A(1)=\emptyset$.
 Suppose $u,v\in A(3)$, $N_{A}(u)=\{u_{1},u_{2},u_{3}\}$ and $N_{A}(v)=\{v_{1},v_{2},v_{3}\}$.
  Then $\dist_{A}(u_{i},u_{j})$ and $\dist_{A}(v_{i},v_{j})$ must be even; otherwise $G$ is not bipartite.
  So, if $N_{A}(v)\cap N_{A}(u)\neq \emptyset$, then $N_{A}(u)=N_{A}(v)$, which implies $|A(3)|\leq 2$.
   If $|A(3)|=1$, then $G \lhd G_{3}$ of Fig. \ref{pic-rank-6}.
   If $|A(3)|=2$, let $u,v\in A(3)$.
   If $N_A(u)=N_A(v)$, then $u \sim v$ as $G$ is reduced by Lemma \ref{neigh-nonaj}; but in this case $G$ would have triangles.
   So $N_A(u) \cap N_A(v)=\emptyset$.
   Noting that $r(G_{3})=6$, $u \nsim v$ in $G$ by Lemma \ref{op1}, which implies that $G= G_{3}$, where $G_3$ is listed in Fig. \ref{pic-rank-6}.

{\it Case 3:}  $ A(1)\neq \emptyset$ and $A(3)\neq \emptyset$.
Assume $u\in A(3)$ and $v\in A(1)$. Noting that $r(A_{4})=r(A_{5})=8$ while $r(G_{1})=6$, we have $\dist_{G}(u,v)=2$, where $A_4,A_5$ are listed in Fig. \ref{picA}.
Then $G \lhd G_{1}$.\hfill$\blacksquare$

\begin{lemma}  \label{bi-B}
Let $G$ be a connected reduced bipartite graph of rank $6$, which contains $B$ of Fig. \ref{pic-nons} as an induced subgraph.
Then $G$ is an induced subgraph of $G_{1}$, $G_{3}$, $G_{4}$ or $G_{5}$.
\end{lemma}

\noindent{\bf Proof.} By Lemma \ref{dist} and Lemma \ref{degree}, $1\leq |N_{B}(v)|\leq 2$ for $v\in V(G)\backslash V(B)$.

{\it Case 1:}  $B(2)\neq \emptyset$.
Let $v\in B(2)$ and let $v_{1}$ and $v_{2}$ be two neighbors of $v$ in $B$.
As $G$ is bipartite, $\dist_{B}(v_{1},v_{2})$ should be even, which implies $\dist_{B}(v_{1},v_{2})=2$.
As $G$ is reduced, $G \rhd A $ (or $\C_6$).
By Lemma \ref{bi-A}, $G$ is an induced subgraph of $G_{1}$, $G_{2}$ or $G_{3}$.
Noting that $B$ is not an induced subgraph of $G_{2}$, so $G\lhd G_{1}$ or $G\lhd G_{3}$.

\begin{figure}[h!]
 \renewcommand\thefigure{\arabic{section}.\arabic{figure}}
  \centering
  \includegraphics[scale=0.7]{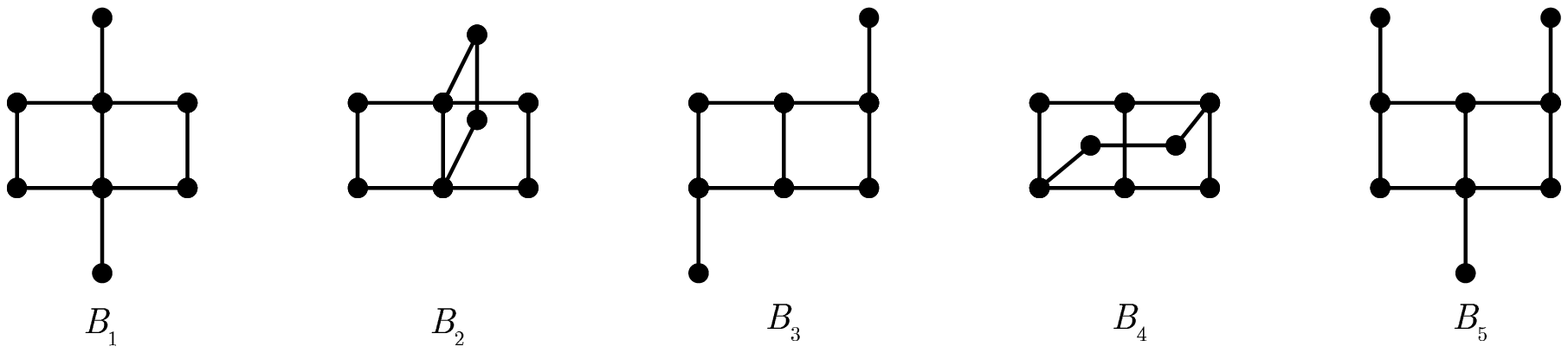}\\
  \caption{\small Illustration for Lemma \ref{bi-B}, where $r(B_{1})=r(B_{2})=r(B_{3})=r(B_{4})=8; r(B_{5})=6$}\label{picB}
\end{figure}

{\it Case 2:} $ B(2)=\emptyset$.
 One can check that $r(B_{i})=8$ for $i=1,2,3,4$, where $B_i$'s are listed in Fig. \ref{picB},
 which shows $|B(1)| \le 3$ by the rank of graphs $B_1,B_3$. If $|B(1)|=1$, the result holds obviously.

 If $|B(1)|=2$, let $u,v \in B(1)$ and let $u',v'$ be their unique neighbors in $B$ respectively.
 As $r(B_3)=r(B_4)=8$, $\dist_{B}(u',v') \le 2$. We have four subcases:

(a) $\dist_{B}(u',v')=1$ and $|N_{B}(u')|=|N_{B}(v')|=2$. As $r(G_{5})=6$, $u \nsim v$ by Lemma \ref{op1}. Then $G\lhd G_{5}$.

(b) $\dist_{B}(u',v')=1$ and one of $u',v'$ has $3$ neighbors in $B$. As $r(B_1)=r(B_2)=8$, exactly one of $u',v'$ has $3$ neighbors in $B$.
As $r(G_{5})=6$, $u\sim v$ by Lemma \ref{op1}. Then $G\lhd G_{5}$.

(c) $\dist_{B}(u',v')=2$ and $|N_{B}(u')|=|N_{B}(v')|=2$. As $r(B_{5})=6$, $u \nsim v$. Then $G\lhd B_{5}\lhd G_{1}$.

(d) $\dist_{B}(u',v')=2$ and exactly one of $\{u',v'\}$ has $3$ neighbors in $B$. As $r(G_{5})=6$, $u \nsim v$. Then $G\lhd G_{5}$.

If $|B(1)|=3$, let $u,v,w \in B(1)$ and $u',v',w'$ be their unique neighbors in $B$ respectively. There are three possible subcases:

(a) $u'\sim v'$ and $v' \sim w'$, $|N_{B}(v')|=3$ and $|N_{B}(u')|=|N_{B}(w')|=2$. As $r(G_{4})=6$, $G\lhd G_{4}$.

(b) $u'\sim v'$ and $v' \sim w'$, $|N_{B}(u')|=3$ and $|N_{B}(v')|=|N_{B}(w')|=2$. As $r(G_{5})=6$, $G\lhd G_{5}$.

(c) $G=B_5$. As $r(B_{5})=6$, $G\lhd G_{1}$.\hfill$\blacksquare$

\begin{lemma}\label{bi-C}
Let $G$ be a connected  reduced bipartite graph with rank $6$, which contains $C$ of Fig. \ref{pic-nons} as an induced subgraph.
 Then $G$ is an induced subgraph of $G_{1}$, $G_{3}$, $G_{4}$ or $G_{5}$.
 \end{lemma}

\noindent{\bf Proof.} We have $1\leq |N_{C}(v)|\leq 2$ for each $v\in V(G)\backslash V(C)$.

{\it Case 1:}  $C(2)\neq \emptyset$.
Let $v\in C(2)$, and let $v_{1},v_{2}$ be its two neighbors in $C$.
As $G$ is bipartite,  $\dist_{C}(v_{1},v_{2})$  should be even,  and $\dist_{C}(v_{1},v_{2})=2$.
As $G$ reduced, one of $\{v_{1},v_{2}\}$ is a pendant vertex of $C$. Thus $G\rhd B$ and the result follows by Lemma \ref{bi-B}.

\begin{figure}[h!]
  \renewcommand\thefigure{\arabic{section}.\arabic{figure}}
  \centering
  \includegraphics[scale=0.7]{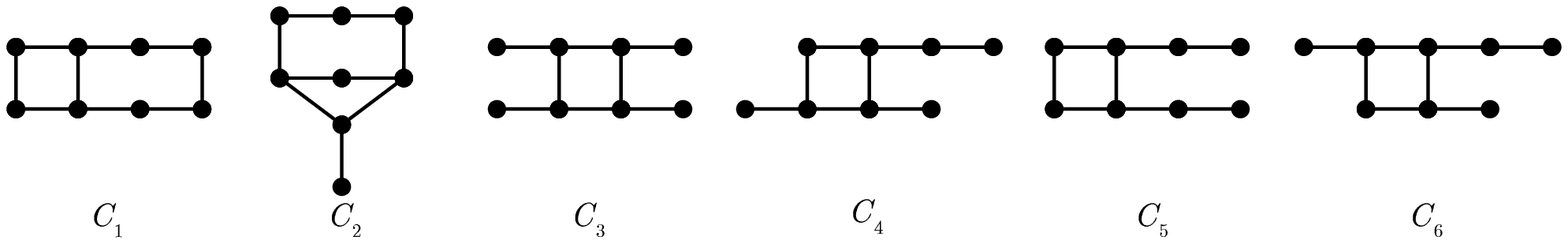}\\
  \caption{\small Illustration for Lemma \ref{bi-C}, where $r(C_{1})=r(C_{2})=r(C_{3})=r(C_{4})=r(C_{5})=8;r(C_{6})=6$}\label{picC}
\end{figure}

{\it Case 2:} $C(2)=\emptyset$.
If $|C(1)|=1$, the result holds obviously. Now we suppose $|C(1)|\geq 2$.
If there exist $u,v\in C(1)$ and $u\sim v$. Let $u'$, $v'$ be their unique neighbors in $C$ respectively.
Then $\dist_{C}(u',v')$ should be odd by the bipartiteness of $G$.
Note that  $|N_{C}(v')|\neq 3$ and $|N_{C}(u')|\neq 3$; otherwise $G$ is not reduced.
If $\dist_{C}(u',v')=1$, then $G\rhd B$.  If $\dist_{C}(v',u')=3$, noting that $r(C_{1})=r(C_{2})=8$, this case cannot occur, where $C_1,C_2$ are listed in Fig. \ref{picC}.

Now assume $v_{1},v_{2},\ldots,v_{k}\in C(1)$ and $v_{i} \nsim v_{j}$ for all $1\leq i< j\leq k$. Similarly, none of them has a neighbor with degree $3$ in $C$. Note that $r(C_{3})=r(C_{4})=8$ and $r(C_{5})=r(C_{6})=6$. We conclude that $G\lhd C_{5}\lhd G_{5}$ or $G\lhd C_{6}\lhd G_{1}$. \hfill$\blacksquare$

\begin{lemma} \label{bi-D}
Let $G$ be a connected reduced bipartite graph with rank $6$, which contains $D$  of Fig. \ref{pic-nons} as an induced subgraph.
 Then $G$ is an induced subgraph of $G_{1}$, $G_{3}$, $G_{4}$ or $G_{5}$.
 \end{lemma}

\noindent{\bf Proof.} We have $1\leq |N_{C}(v)|\leq 3$ for each $v\in V(G)\backslash V(C)$.

{\it Case 1:} $D(2)\neq \emptyset$.
 Let $v\in D(2)$, and let $v_{1}$ and $v_{2}$ be its two neighbors in $D$.
 The $\dist_{D}(v_{1},v_{2})$ should be even.
 If $\dist_{D}(v_{1},v_{2})=2$, then $G$ is not reduced.
 If $\dist_{D}(v_{1},v_{2})=4$, then $G\rhd A$. Thus the result follows from Lemma \ref{bi-A}.

{\it Case 2:} $D(3)\neq \emptyset$.
Let $v\in B(3)$, and let $v_{1},v_{2}$ and $v_{3}$ be its three neighbors in $D$.
Since $\dist_{D}(v_{i},v_{j})$ must be even, we have $G\rhd B$. Then the result follows from Lemma \ref{bi-B}.

{\it Case 3:} $D(2)=D(3)=\emptyset$.
Let $v_{1},v_{2},\ldots,v_{k}\in D(1)$ and let $v_{1}',v_{2}'\ldots,v_{k}'$ be their unique neighbors in $D$ respectively.
As $G$ is reduced, $v_{i}'$ cannot be the quasi-pendant vertex of $D$ (i.e., the vertex adjacent to the pendant vertex).
If $k=1$, then $G\lhd G_{1}$ and the result follows easily.

If $v_{i}\sim v_{j}$ for some $1\leq i\leq j\leq k$, then $\dist_{D}(v_{i}',v_{j}')$ should be odd. We have the following  subcases:


(a) If $\dist_{D}(v_{i}',v_{j}')=1$, then $v_{i}'$ and $v_{j}'$ must be the two vertices in the middle of the path $D$. Thus $G\rhd C$ and the result follows.

(b) If $\dist_{D}(v_{i}',v_{j}')=3$, then one of $\{v_{i}',v_{j}'\}$ must be the pendant vertex of $D$. Thus $G\rhd A$ and the result follows.

\begin{figure}[h!]
 \renewcommand\thefigure{\arabic{section}.\arabic{figure}}
  \centering
  \includegraphics[scale=0.7]{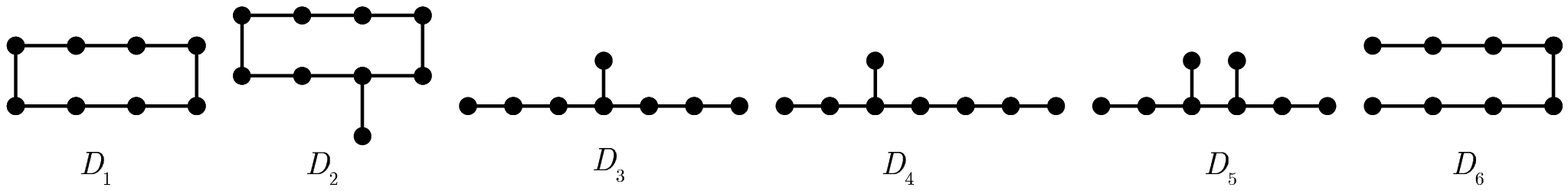}\\
  \caption{\small Illustration for Lemma \ref{bi-D}, where $r(D_{1})=r(D_{3})=6;r(D_{2})=r(D_{4})=r(D_{5})=r(D_{6})=8$}\label{picD}
\end{figure}

(c) If $\dist_{D}(v_{i}',v_{j}')=5$, then $G\rhd D_{1}$. Noting that $r(D_{1})=6$ and $r(D_{2})=8$, we have $G\lhd D_{1}\lhd G_{4}$, where $D_1,D_2$ are listed in Fig. \ref{picD}.

 If $v_{i} \nsim v_{j}$ for any $1\leq i< j\leq k$. Noting that $r(D_{3})=6$ and $r(D_{4})=r(D_{5})=r(D_{6})=8$,  We have $G\lhd D_{3}\lhd G_{1}$,
 where $D_3,D_4,D_5,D_6$ are listed in Fig. \ref{picD}. \hfill$\blacksquare$

\begin{lemma} \label{bi-E}
Let $G$ be a connected reduced bipartite graph with rank $6$, which contains $E$ of Fig. \ref{pic-nons} as an induced subgraph.
Then $G$ is an induced subgraph of $G_{1}$, $G_{2}$ $G_{3}$, $G_{4}$ or $G_{5}$.
\end{lemma}

\noindent{\bf Proof.} We have $1\leq |N_{E}(v)|\leq 3$ for each $v\in V(G)\backslash V(E)$.

{\it Case 1:} $E(2)\neq \emptyset$.
Let $v\in E(2)$ and let $v_{1}$ and $v_{2}$ be its two neighbors in $E$.
The $\dist_{E}(v_{1},v_{2})$ is even.
If $\dist_{E}(v_{1},v_{2})=4$, then $G\rhd A$.
 Suppose $\dist_{E}(v_{1},v_{2})=2$.
Noting that $G$ is reduced, $G\rhd A$ or $G\rhd C$. The result follows from Lemmas \ref{bi-A} and \ref{bi-C}.

{\it Case 2:} $E(3)\neq \emptyset$.
 Let $v\in E(3)$ and let $v_{1},v_{2}$ and $v_{3}$ be the neighbors of $v$ in $E$. The $\dist_{E}(v_{i},v_{j})$ is even for $1\leq i< j\leq 3$.
 By Lemma \ref{neigh-outside}, $G\rhd B$ and the result follows.

\begin{figure}[h!]
  \renewcommand\thefigure{\arabic{section}.\arabic{figure}}
  \centering
  \includegraphics[scale=0.7]{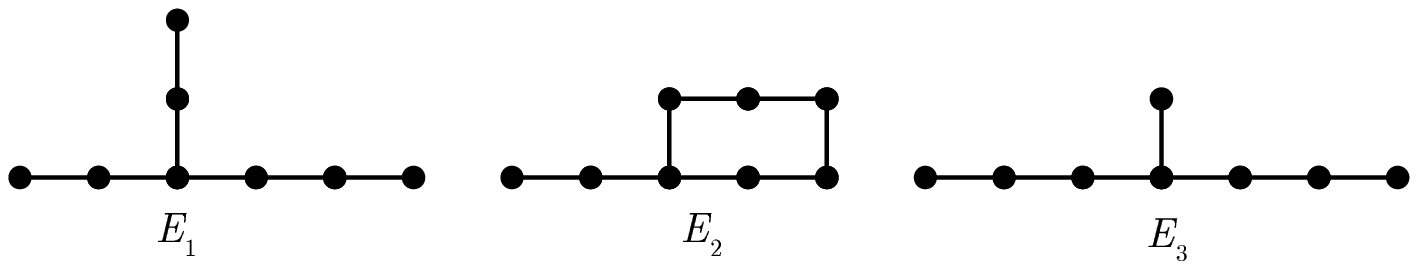}\\
  \caption{\small Illustration for Lemma \ref{bi-E}, where $r(E_{1})=r(E_{2})=8;r(E_{3})=6$}\label{picE}
\end{figure}

{\it Case 3:} $E(2)=E(3)=\emptyset$.
For each $v\in E(1)$, it must be adjacent to a pendant vertex of $E$ since $G$ is reduced,
which implies $|E(1)|\leq 3$. Let $v\in E(1)$,  and let $v'$ be its unique neighbor in $E$ (i.e. a pendant vertex of $E$).
If $|E(1)|=1$, then $G\lhd G_{1}$ and the result follows.
If $|E(1)|\geq 2$, noting that $r(E_{1})=r(E_{2})=8$ and $r(E_{3})=6$, we have $G\lhd E_{3}\lhd G_{1}$, where $E_1,E_2,E_3$ are listed in Fig. \ref{picE}. \hfill$\blacksquare$

\begin{lemma} \label{p4p2}
Let $G$ be a connected reduced bipartite graph with rank $6$, which contains $\P_{4}\cup \P_{2}$  of Fig. \ref{pic-nons} as an induced subgraph.
Then $G$ is an induced subgraph of $G_{1},G_{2},G_{4}$ or $G_{5}$.
\end{lemma}

\noindent{\bf Proof.}
Note that $G$ is connected and bipartite, and $\dist_{G}(v, \P_{4}\cup \P_{2})=1$ for each vertex $v\in V(G)\backslash V(\P_{4}\cup \P_{2})$. We have two cases:

{\it Case 1:}  There exists a vertex $v\in V(G)\backslash V(\P_{4}\cup \P_{2})$ such that $\dist_{G}(v,\P_{4})=\dist_{G}(v,\P_{2})=1$.
Then $|N_{\P_{2}}(v)|=1$ and $|N_{\P_{4}}(v)|\leq 2$. If $|N_{\P_{4}}(v)|=1$, then $G\rhd D$; if $|N_{\P_{4}}(v)|=2$, the $G\rhd C$.
 The result follows by Lemmas \ref{bi-D} and \ref{bi-C}.

{\it Case 2:} There exist two vertices $u,v\in V(G)\backslash V(\P_{4}\cup \P_{2})$ such that $\dist_{G}(u,\P_{4})=\dist_{G}(v,\P_{2})=1$ and $u\sim v$.
Then $G\rhd D$ and the result also follows.\hfill$\blacksquare$

\begin{lemma} \label{3p2}
Let $G$ be a connected reduced bipartite graph with rank $6$, which contains $3\P_{2}$ of Fig. \ref{pic-nons} as an induced subgraph.
Then $G$ is an induced subgraph of $G_{1},G_{2},G_{3},G_{4}$ or $G_{5}$.
\end{lemma}

\noindent{\bf Proof.} We denote the three disjoint paths as $\P_{2}^{1},\P_{2}^{2},\P_{2}^{3}$ respectively. Then we have two cases:

{\it Case 1:}  There exist two vertices $u,v\in V(G)\backslash V(3\P_{2})$ such that $\dist_{G}(u,\P_{2}^{1})=\dist_{G}(v,\P_{2}^{2} )=1$ and $u\sim v$.
Then $G\rhd D$ and the result follows by Lemma \ref{bi-D}.

{\it Case 2:} There exists a vertex $v\in V(G)\backslash V(3\P_{2})$ such that $\dist_{G}(v,\P_{2}^{1})=\dist_{G}(v,\P_{2}^{2})=1$.
 If $\dist_{G}(v,\P_{2}^{3})=1$, then $G\rhd E$;
 if $\dist_{G}(v,\P_{2}^{3})\geq 2$, then $G\rhd (\P_{4}\cup \P_{2})$.
 The result follows by Lemmas \ref{bi-E} and \ref{p4p2} .\hfill$\blacksquare$

\subsection{Reduced triangle-free and non-bipartite graphs with rank $6$}
Observe that if a graph $G$ contains an induced odd cycle $\C_{2k+1}$, then $r(G) \ge r(\C_{2k+1})=2k+1$.
So, if $G$ is triangle-free, non-bipartite, and has rank $6$, then $G$ contains an induced $\C_5$.

\begin{lemma} \label{nonbi-rank6}
Let $G$ be a triangle-free and non-bipartite graph with rank  $6$.
Then $G\rhd F$ of Fig. \ref{pic-nons}.\end{lemma}

\noindent{\bf Proof.}
As discussed above, $G$ contains an induced cycle $\C_5$.
Let $v$ be an arbitrary vertex not in $\C_{5}$.
Note that $r(\C_{5})=5$, thus $d_{\C_{5}}(v)\geq 1$ by Lemma $\ref{dist}$.
If $d_{\C_{5}}(v)=1$, then $G\rhd F$. If $d_{\C_{5}}(v)\geq 3$, then $G$ contains triangles; a contradicition.
Now suppose $d_{\C_{5}}(v)=2$ for each $v\notin V(\C_{5})$.
Let $V(G)\backslash V(\C_{5})=\{v_{1},v_{2},\ldots,v_{k}\}$.
If $v_{i_{0}}\sim v_{j_{0}}$ for some $i_{0},j_{0}$, then $N_{\C_{5}}(v_{i_{0}})\cap N_{\C_{5}}(v_{j_{0}})=\emptyset$; otherwise $G$ contains triangles.
 Thus for any $i \ne j$, either $v_{i}$ is not adjacent to $v_{j}$, or $v_{i}\sim v_{j}$ and $N_{\C_{5}}(v_{i})\cap N_{\C_{5}}(v_{j})=\emptyset$.
 However, in any case $G$ can be obtained from $\C_{5}$ by several steps of multiplication of vertices, which implies $r(G)=r(\C_{5})=5$.\hfill$\blacksquare$

\begin{figure}[h!]
  \renewcommand\thefigure{\arabic{section}.\arabic{figure}}
  \centering
  \includegraphics[scale=0.7]{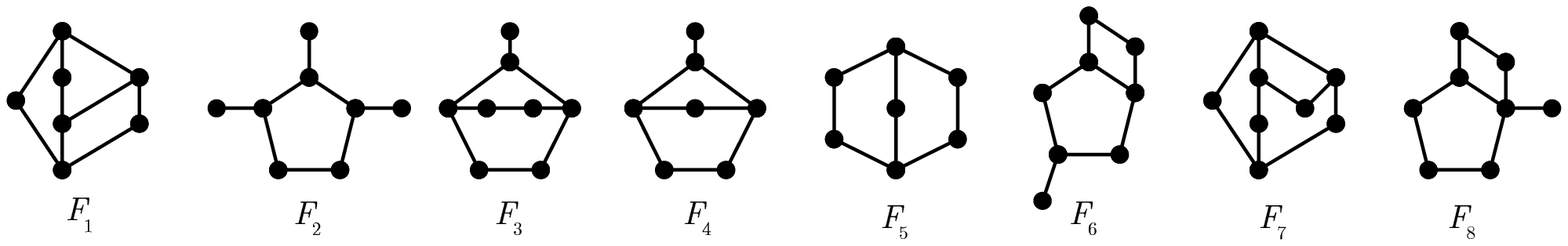}\\
  \caption{\small Illustration for Lemma \ref{red-nonbi-rank6}, where $r(F_{i})\geq 7$, $i=1,2,\ldots,8$}\label{picF}
\end{figure}

\begin{lemma} \label{red-nonbi-rank6}
Let $G$ be a connected reduced triangle-free graph with rank $6$, which contains $F$ of Fig. \ref{pic-nons} as an induced subgraph.
 Then $G$ is an induced subgraph of $G_{6}$ of Fig. \ref{pic-rank-6}.
 \end{lemma}

\noindent{\bf Proof.} If $|N_{F}(v)|\geq 4$ for some $v\in V(G)\backslash V(F)$, then $G\rhd C_{3}$.
If $|N_{F}(v)|=3$ for some $v\in V(G)\backslash V(F)$, say $N_{F}(v)=\{v_{1},v_{2},v_{3}\}$,
then $\dist_{F}(v_{i},v_{j})\geq 2$ for $1\leq i<j\leq 3$.
Thus one of $\{v_{1},v_{2},v_{3}\}$ is the pendant vertex in $F$.
So, $G\rhd F_{1}$ or $G$ is obtained from $F$ by multiplication of the vertex with degree $3$.
However, $r(F_{1})=7$, where $F_1$ is listed in Fig. \ref{picF}.
Now we conclude that $1\leq |N_{F}(v)|\leq 2$ for $v\in V(G)\backslash V(F)$.

{\it Case 1:} $F(2)=\emptyset$.
Let $F(1)=\{v_{1},v_{2},\ldots,v_{k}\}$, and $N_{F}(v_{i})=\{u_{i}\}$ for $i=1,2,\ldots,k$.
Then $u_{i}$ cannot be adjacent to the pendant or the quasi-pendant vertex of $F$;
otherwise $r(G)\geq 7$ by Lemma \ref{pend-nul} or $G$ is not reduced by Lemma \ref{neigh-outside}.
Thus $|F(1)|\leq 4$. One can check that $r(F_{2})=r(F_{3})=8$, and $r(G_{6})=6$. By Corollary \ref{op2} we have $G\lhd G_{6}$.

{\it Case 2:}  $F(1)=\emptyset$.
Let $w_{1}$ be the vertex with maximum degree in $F$, and $N_{F}(w_{1})=\{w_{2},w_{3},w_{4}\}$, where $w_{2}$ is the pendant vertex of $F$.
Let $F(2)=\{v_{1},v_{2},\ldots,v_{k}\}$ and $N_{F}(v_{i})=\{u_{i}^{1},u_{i}^{2}\}$ for $i=1,2,\ldots,k$.
Since $G$ is triangle-free, we have $2\leq \dist_{F}(u_{i}^{1},u_{i}^{2})\leq 3$.
If $u_{i}^{1}=w_{3}$ and $u_{i}^{2}=w_{4}$ for some $i$, then $r(G)\geq r(F_{4})=7$, where $F_4$ is listed in Fig. \ref{picF}.
As $G$ is reduced,  one of $\{u_{i}^{1},u_{i}^{2}\}$ must be $w_{2}$.
Set $F(2)=F(2)_{1}\cup F(2)_{2}$, where $F(2)_{1}=\{v_{i}|\dist_{F}(u_{i}^{1},u_{i}^{2})=2\}$ and $F(2)_{2}=\{v_{i}|\dist_{F}(u_{i}^{1},u_{i}^{2})=3\}$.
As $r(F_{5})=7$, then $F(2)_{2}=\emptyset$, where $F_5$ is listed in Fig. \ref{picF}.
 Note that $|F(2)_{1}|\leq 2$.
 Suppose that $F(2)_{1}=\{v_{1},v_{2}\}$.
 We have $N_{F}(v_{1})=\{w_{2},w_{3}\}$ and $N_{F}(v_{2})=\{w_{2},w_{4}\}$.
 If $v_{1}\sim v_{2}$, $G\rhd C_{3}$; otherwise $G\rhd C_{7}$ which implies $r(G)\geq 7$.
 Thus $|F(2)_{1}|\leq 1$ and $G\lhd G_{6}$.

{\it Case 3:}   $ F(1)\neq\emptyset$ and $F(2)\neq\emptyset$.
By the above discussion, we know $|F(2)_{1}|=1$ and $F(2)_{2}=\emptyset$.
By a simple checking, $r(F_{i})=8$ for $i=6,7,8$.
So $G\lhd G_{6}$.\hfill$\blacksquare$

\end{document}